%% file: lamhypbisarxiv.tex
\documentclass[12pt,a4paper,leqno]{amsart} 
\ifx\pdftexversion\undefined
  \usepackage[dvips]{graphicx,color}
\else
  \usepackage[pdftex]{graphicx,color}
\fi
\usepackage[francais]{babel}
\usepackage[latin1]{inputenc}%

\let\ifanglais\iffalse

\def\R{{\mathbb R}}

\def\vole{{\rm vol}_e}
\def\N{{\mathbb N}}

\newtheoremstyle{mesthm}
  {10pt plus 1pt minus 1pt}
  {9pt minus 6pt}
{\slshape}
  {0.5cm}
  {\bfseries}
  {.}
  {1ex}
  {}
\newtheoremstyle{mesdefi}
  {6pt plus 1pt minus 1pt}
  {6pt plus 1pt minus 1pt}
  {}
  {0.5cm}
  {\bfseries}
  {.}
  {1ex}
  {}%

\theoremstyle{mesthm}
\newtheorem{lemm}{\ifanglais{\large L}emma\else{\large L}emme\fi}[section]
\newtheorem{theo}[lemm]{\ifanglais{\large T}heorem\else {\large
    T}héorème\fi}

\newtheorem{prop}[lemm]{{\large P}roposition} 
 
\newtheorem{coro}[lemm]{\ifanglais{\large C}orollary\else{\large C}orollaire\fi}

\theoremstyle{mesdefi}
\newtheorem{defi}[lemm]{\ifanglais{\large D}efinition\else{\large
    D}éfinition\fi} 
 \newtheorem{rema}[lemm]{\ifanglais{\large
    R}emark\else{\large  R}emarque\fi}

\makeatletter
\renewenvironment{proof}[1][{\sl \ifanglais Proof\else Démonstration\fi}]{\par
  \normalfont
  \topsep6\p@\@plus6\p@ \trivlist
  \item[\hskip\labelsep\slshape
    #1\@addpunct{.}]\ignorespaces
}{%
  \qed\endtrivlist
}
\makeatother

\newcommand{\n}[1]{\Vert #1 \Vert}

\title[Géométrie locale des géométries de Hilbert]%
{Les géométries de Hilbert sont à géométrie locale bornée}

\author[B.~Colbois et C.~Vernicos]{Bruno Colbois et Constantin Vernicos}
\address{Institut de mathématique\\
  Université de Neuchâtel\\
  Rue Émile Argand 11\\
  Case postale 158\\
  CH--2009 Neuchâtel\\
  Switzerland}
\email{Bruno.Colbois@unine.ch\\
Constantin.vernicos@unine.ch}

\subjclass[2000]{Differential Geometry; Metric Geometry}

\begin{document}

\maketitle

\begin{abstract}
  On montre que la géométrie de Hilbert d'un domaine convexe de $\R^n$ est
  à géométrie locale bornée c-à-d que pour un rayon fixé, toutes les boules
sont bilipschitz à un domaine de $\R^n$ euclidien. On en déduit que si
la géométrie de Hilbert est hyperbolique au sens de Gromov, alors le bas de
son spectre est strictement positif. On donne un contre-exemple en dimension trois qui montre
que la réciproque n'est pas vraie pour les géométries de Hilbert non planes.
\end{abstract}

\renewcommand{\abstractname}{Abstract}
\begin{abstract}(\textsl{Hilbert geometries have bounded local geometry})
We prove that the Hilbert geometry of a convex domain in $\R^n$ has bounded local geometry,
i.e., for a given radius, all balls are bilipschitz to a euclidean domain of $\R^n$.
As a consequence, if the Hilbert geometry is also 
Gromov hyperbolic, then the bottom of its spectrum is strictly positive. We also
give a counter exemple in dimension three wich shows that the reciprocal is not true
for non plane Hilbert geometries.
\end{abstract}

\input{spectre05}

\nocite{*}

\bibliographystyle{amsalpha}
\bibliography{airhyp}

\end{document}

%% file: spectre05.tex

\section{Introduction}

Le but de cet article est de montrer que localement, toute géométrie de Hilbert
est proche de la géométrie euclidienne, et cela de manière uniforme. En géométrie métrique,
le contrôle de la géométrie locale joue un rôle croissant dans l'étude
des quasi-isométries (voir par exemple le concept
«\,d'espace métrique à croissance bornée à une certaine échelle\,», \cite{bs}~p.~294, ou encore \cite{lp}, paragraphe 5). En 
géométrie riemannienne, il est bien connu que  la
géométrie locale (notamment par l'intermédiaire de la courbure et du rayon d'injectivité) 
intervient dans la démonstration
 de nombreux théorèmes. Cependant, si les géométries de Hilbert sont 
naturellement munies d'une métrique de Finsler, celle-ci est seulement $C^0$,
ce qui exclu tout concept de courbure proche de la notion riemannienne 
correspondante. Aussi, nous controlerons localement la géométrie d'un point de vue métrique
en demandant que les boules de rayon fixé soient uniformément bilipschitz
 à une partie de l'espace euclidien. Nous en déduirons que l'hyperbolicité
au sens de Gromov d'une géométrie de Hilbert implique la non-nullité du 
bas du spectre: cela avait été démontré en dimension 2 dans \cite{cv1}, 
et était conjecturé en dimension supérieure. La preuve que nous donnons ici
est nouvelle également en dimension 2, et plus simple que la précédente. 
Pour conclure avec cette question, 
nous donnons un exemple montrant que contrairement à la dimension 2, il n'y a
pas de
réciproque à ce dernier résultat en dimension 
supérieure.

Avant d'énoncer les résultats précis, rappelons qu'une géométrie de Hilbert
$(\mathcal C, d_{\mathcal C})$ est la donnée d'un ouvert convexe et borné 
$\mathcal C$ de $\mathbb R^n$ muni de la distance de Hilbert $d_{\mathcal C}$ définie
de la manière suivante : 
pour toute paire de points distincts $p$ et $q$ dans $\mathcal{C}$, la droite passant
par $p$ et $q$ rencontre le bord $\partial \mathcal C$ de 
$\mathcal C$ en deux points distincts $a$ et $b$ tels que la droite passe
par $a$, $p$, $q$ et $b$ dans cet ordre.

On définit alors
$$d_{\mathcal C}(p,q) = \frac{1}{2} \ln [a,p,q,b]$$
où $[a,p,q,b]$ est le birapport de $(a,p,q,b)$:
$$[a,p,q,b] = \frac{\Vert q-a \Vert_{e} }{\Vert p-a \Vert_{e}} × 
\frac{\Vert 
p-b \Vert_{e}}{\Vert q-b \Vert_{e}}   $$
où $ \Vert . \Vert_{e}$ désigne la norme euclidienne.
On pose également $d_{\mathcal C}(p,p) =0$.

Sur $\mathcal C$, on peut mettre une norme de Finsler $C^0$, 
notée  $\n{\cdot}_{\mathcal C}$, en procédant
comme suit: si $p \in \mathcal C$ et $u_{p} \in T_{p}\mathcal C = \mathbb 
R^n$, la droite passant par $p$ et dirigée par $u_{p}$ coupe $\partial \mathcal 
C$ en deux points $p^{+}$ et $p^-$.  
On pose alors
$$\n{u_{p}}_{\mathcal C} = \frac{1}{2} \n{u_{p}}_{e} 
\biggl(\frac{1}{\Vert 
p-p^{+} \Vert_{e}} + 
\frac{1}{\Vert p-p^- \Vert_{e}}\biggr) {\rm }$$

où $ \n{u_{p}}_{e}$ désigne la norme euclidienne de $u_{p}$.

De plus, on utilisera dans la suite le fait que la distance de longueur induite sur $\mathcal C$ par la norme 
$\Vert . \Vert_{\mathcal C}$ coïncide avec $d_{\mathcal C}$ --- voir \cite{so0} pour une justification
détaillée de cette dernière affirmation.

A cette norme de Finsler est associée une norme duale: si $l_{p}$ 
est une forme linéaire sur $T_{p}\mathcal C$, on pose
$$\n{l_{p}}_{\mathcal C}^{*} = \sup \{l_{p}(u_{p}): 
u_{p}\in T_{p}\mathcal C,\ \n{u_{p}}_{\mathcal C}=1 \}.$$

Grâce à la norme de Finsler, on construit une forme volume et 
une mesure sur 
$\mathcal C$ (qui correspond en fait à  
la mesure de Hausdorff, voir \cite{bbi} 
exemple 5.5.13).

Soient 
$p \in \mathcal C$ et 
$$TB_{\mathcal C}(p) = \{u_{p} \in T_{p}\mathcal 
C =\mathbb R^n : 
\n{u_{p}}_{\mathcal C} < 1 \}$$
la boule unité de $T_{p}\mathcal C$. 
Soit $\omega_{n}$ le volume de la boule unité de l'espace euclidien
$\mathbb R^n$ 
et considérons la fonction $h_\mathcal{C}\colon\mathcal C \longrightarrow \mathbb R$ donnée par 
$$
{h_\mathcal{C}(p)=\omega_{n}/\vole\bigl(TB_{\mathcal C}(p)\bigr)}
$$
où $\vole$ est le volume euclidien usuel. 

Alors la mesure $\mu_{\mathcal C}$ (que nous nommerons   
mesure de Hilbert de $\mathcal{C}$; elle est également connue sous le nom de mesure
de Busemann)  associée à $\n{\cdot}_{\mathcal C}$ 
est définie ainsi : si $A \subset \mathcal C$ est un borélien, on pose
$$\mu_{\mathcal C} (A) = \int_{A} h_\mathcal{C}(p)d\vole(p)$$
où $d\vole(p)$ est la mesure de Lebesgue.

Si $\Omega$ est un domaine avec $\overline{\Omega} \subset\mathcal C$, l'intégrale 
sur $\Omega$ par 
rapport à $\mu_{\mathcal C}$ d'une 
fonction $f$ définie sur $\Omega$ sera notée
$\int_\Omega f~d\mu_\mathcal{C}.$

\begin{rema}
  Les résultats de cet article restent inchangés si on considère une mesure
équivalente. En particulier, la mesure dite de Holmes-Thompson, qui
est équivalente par les inégalités de Santalo et Bourgain-Milman à la mesure de Hilbert, donne
les mêmes résultats. Remarquons que bien que cette mesure soit
liée à un invariant symplectique, c'est une mesure à densité par
rapport à la mesure de Lebesgue, et l'expression de cette
densité peut toujours être définie dans le cadre des géométries
de Hilbert.
\end{rema}

On définit le bas du spectre de $\mathcal C$, que l'on note 
$\lambda_1(\mathcal C)$, par analogie avec ce qui se fait dans le cas 
des variétés riemanniennes de volume infini. On pose
\begin{equation}
\label{eqdef1}
\lambda_1(\mathcal{C})=
\inf \frac{\displaystyle{\int_\mathcal{C}{\n{df_p}_\mathcal{C}^*}^2~d\mu_\mathcal{C}(p)}}
{\displaystyle{\int_\mathcal{C} f^2(p)d\mu_\mathcal{C}(p)}}\text{,}
\end{equation}
où l'infimum est pris sur toutes les fonctions 
lipschitziennes, non nulles, à support compact dans $\mathcal{C}$ et où
$\mu_\mathcal{C}$ est la mesure de Hilbert associée à $\mathcal C$. 
 L'expression ci-dessus est appelée le quotient de Rayleigh de $f$.
 
 La constante de Cheeger de $\mathcal C$ est définie par
 \begin{equation}
\label{Cheeger}
I^{\infty}(\mathcal{C})=
\inf_U \frac{\bar{\nu}_{\mathcal C}(\partial U)}{\mu_{\mathcal C}(U)}\text{,}
\end{equation}
où $U$ décrit les ouverts de $\mathcal C$ d'adhérence compacte dont le 
bord $\partial U$ est une sous-variété de dimension n-1. On rappelle que
$\bar{\nu}_{\mathcal C}$ désigne ici la mesure de Hausdorff des hypersurfaces
associée à la restriction de $\Vert . \Vert_{\mathcal C}$.

Remarquons qu'aussi bien le bas du spectre que la constante de Cheeger ne dépendent
que $C^0$ de la métrique.

Notre résultat principal est

\begin{theo}\label{theoprincipal}
  Un convexe $(\mathcal{C},d_\mathcal{C})$ muni de sa métrique de Hilbert est à géométrie
locale bornée. 
\end{theo}

\begin{defi}\label{geobornee} 
  Nous dirons qu'un espace métrique $(X,d)$ est à \textbf{géométrie locale bornée} s'il
existe deux constantes $a$ et $C$ telles que pour tout $x\in X$, si $B_a(x)$ est la boule centrée en $x$
de rayon $a$, il existe une application
$F\colon B_a(x)\to \R^n$ qui soit un homéomorphisme sur son image et soit 
C-bilipschitz, i.e., pour tout $y$ et $z$ dans $B_a(x)$
$$
\frac{1}{C}d(z,y)\leq \n{F(z)-F(y)}_e \leq C d(z,y)\text{.}
$$ 
\end{defi}

Comme application du théorème \ref{theoprincipal}, on déduit

\begin{theo}\label{equiv1}
Soit $(\mathcal{C},d_\mathcal{C})$  un convexe de $R^n$ muni de sa métrique de Hilbert.
S'il existe $\delta>0$ tel que $(\mathcal{C},d_\mathcal{C})$ est $\delta$-hyperbolique alors son bas
du spectre est non nul, i.e., $\lambda_1(\mathcal{C})> 0$.
\end{theo}

Nous démontrons ce résultat en utilisant la constante de Cheeger et l'inégalité
liant la constante de Cheeger au bas du spectre  démontrée dans \cite{cv1}. Plus précisément nous
démontrons le théorème suivant:

\begin{theo}\label{equiv2}
Soit $(\mathcal{C},d_\mathcal{C})$  un convexe de $R^n$ muni de sa métrique de Hilbert.
S'il existe $\delta>0$ tel que $(\mathcal{C},d_\mathcal{C})$ est $\delta$-hyperbolique alors sa constante
de Cheeger $I^\infty(\mathcal{C})$ est strictement positive.
\end{theo}

A la section 1, on donne la preuve du théorème \ref{theoprincipal}. La section 2
est consacrée aux théorèmes \ref{equiv1} et \ref{equiv2} et, à la section 3, on construit
un exemple montrant que la réciproque aux théorèmes \ref{equiv1} et \ref{equiv2} n'est pas vraie.



\section{Démonstration du théorème \ref{theoprincipal}}

Pour démontrer le théorème \ref{theoprincipal}, on va montrer que la géométrie
des boules de rayon $1$ est bornée (i.e. on peut prendre $a=1$ dans la définition).

L'idée est d'utiliser l'ellipsoïde de John: on sait que les boules d'une géométrie
de Hilbert sont affinement convexes (mais en général pas métriquement)
--- voir \cite{so0} p. 21. On sait qu'à tout convexe borné $C$ de
$\R^n$, on peut associer l'ellipsoïde de John qui est l'unique ellipsoïde $E$
de volume euclidien maximal contenu dans $C$. De plus, on a 
$E \subset C \subset \sqrt n E$. 

\begin{proof}[Démonstration du théorème \ref{theoprincipal}]
  Soit $p$ un point de $\mathcal{C}$ et $B_\mathcal{C}(p,1)$ la boule unité. Considérons 
$\mathcal{B}$ l'ellipsoïde
de John inclu dans $B_\mathcal{C}(p,1)$.

Quitte à translater et à appliquer une transformation affine,
on peut supposer que l'ellipsoïde de John $\mathcal{B}$ est la boule
unité euclidienne centrée en l'origine. 

Ainsi pour tout $y,z \in \partial B_\mathcal{C}(p,1)$ on a
d'une part
$$
1\leq \Vert y \Vert_e \leq \sqrt{n} 
$$
et aussi
$$
d_\mathcal{C}(y,z)\leq 2
$$

Pour contrôler la distance de Hilbert sur $B_\mathcal{C}(p,1)$, on va en fait
encadrer la norme de Finsler $\Vert . \Vert_{\mathcal C}$  et montrer que
$(B_{\mathcal{C}}(p,1), \Vert . \Vert_{\mathcal C})$ est bilipschitz, au sens
de Finsler (i.e., infinitésimalement), à
$B_{\mathcal{C}}(p,1)$ avec la métrique euclidienne usuelle. Cela implique 
immédiatement la même propriété pour $d_{\mathcal C}$ qui est la distance de 
longueur associée.

Pour cela on choisit $x\in B_{\mathcal{C}}(p,1)$ et on estime
$$
\Vert v \Vert_{\mathcal C}=\frac{1}{2}\n{v}_e\left(\frac{1}{\n{x-x_+}_e}+\frac{1}{\n{x-x_-}_e} \right)
$$
où $v$ est un vecteur tangent à $\mathcal C$ au point $x$.

Pour pouvoir comparer à la métrique euclidienne, c'est-à-dire à $\n{v}_e$,
 il faut et suffit 
d'une part que les nombres
$\n{x-x_+}_e$ et $\n{x-x_-}_e$ ne soient pas trop petits et d'autre part 
que l'un d'eux au moins ne soit pas trop grand. Cela se traduit géométriquement par le 
fait que le bord du convexe $\mathcal C$ ne doit pas être trop proche ou trop éloigné
de $B_{\mathcal{C}}(p,1)$ dans un sens qui sera précisé dans la suite.

A l'étape 1, on montre qu'il n'est pas trop proche, à l'étape 2 qu'il
n'est pas trop éloigné, et à l'étape 3, on explique comment cela 
permet de conclure la preuve du théorème. 

\medskip
\textbf{Etape 1:} La distance entre $\partial \mathcal C$ et
$\partial B_{\mathcal{C}}(p,1)$ est minorée par $\frac{1}{2e^4+1}$.

\smallskip
L'idée intuitive est ici que si $\partial \mathcal C$ est trop proche du bord
de la boule unité $B_{\mathcal{C}}(p,1)$, on arrive à exhiber deux points de
la boule qui sont trop éloignés. Dans la preuve ci-dessous, on utilise le fait
que les points de la boule $B_{\mathcal{C}}(p,1)$ ne sont pas distants de plus de $2$
pour en déduire l'estimation explicite.

\smallskip
Soit $q$ un point de $\partial \mathcal{C}$  et $q_0$ un point de $\partial B_\mathcal{C}(p,1)$ 
réalisant le minimum de la distance euclidienne de $\partial \mathcal{C}$
à $\partial B_\mathcal{C}(p,1)$.

Le but est de montrer que $d_0=qq_0$ la distance euclidienne entre ces deux ensembles  
est minorée par $\frac{1}{2e^4+1}$.

On se place dans un plan contenant l'origine et les points $q_0$ et $q$.
Sans limiter la généralité, après une éventuelle rotation ou symétrie, on se trouve
dans la situation de la figure \ref{figdeux}, i.e., $q_0$ est sur l'axe des ordonnées positifs. On
note également $\rho_0$ la distance euclidienne de $q_0$ à $\mathcal{B}$, ainsi les coordonnées de $q_0$
sont $(0,\rho_0+1)$.

\begin{figure}[hbtp]
  \centering
  \includegraphics[scale=.6]{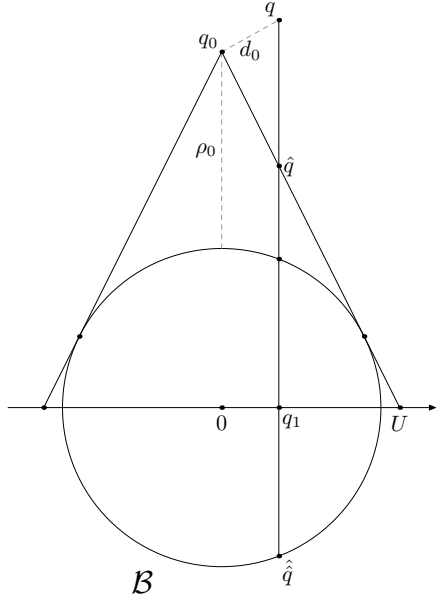}
\caption{Minoration de $d_0$}  \label{figdeux}
\end{figure}

Si $d_0$ est plus grand que $1$, on a la minoration cherchée.

\smallskip
Supposons donc que $d_0<1$. L'idée est que le cône $C_{q_0}$ de sommet $q_0$ et tangent à $\mathcal{B}$
est tout entier dans la boule $ B_\mathcal{C}(p,1)$, par convexité de cette dernière.
Soit $D_q$ la parallèle à $0q_0=0y$, que l'on a pris pour axe des ordonnées, passant par $q=(q_1,q_2)$. 
Puisque $d_0<1$, la droite $D_q$ coupe la boule $\mathcal{B}$
en deux points. On note $\hat{\hat q}=(q_1,-\sqrt{1-q_1^2})$ le point d'ordonnée négative. 
Quant à $\hat q=(q_1,\hat q_2)$, il désignera
le point d'intersection de $D_q$ avec l'une des deux tangentes à $\mathcal{B}$ passant par $q_0$.

Pour minorer $d_0$, on se place sur la droite $D_q$
et on utilise le fait que $ \hat q$ et $\hat{\hat q}$ sont dans $ B_\mathcal{C}(p,1)$ d'où
\begin{equation}
  \label{eqdoro}
  2\geq d_\mathcal{C}(\hat q,\hat{\hat q}) \geq \frac{1}{2} \ln \frac{q_2+\sqrt{1-q_1^2}}{q_2-\hat q_2}\text{.}
\end{equation}
Soit $2U$ la largeur du cône $C_{q_0}$ lorsqu'il rencontre l'axe des abscisses.  
Par le théorème de Thales on a
\begin{equation}\label{eqdoro1}
  \frac{\hat q_2}{1+\rho_0} = \frac{U -q_1}{U}  \text{ donc }  \hat q_2= (1+\rho_0)\frac{U -q_1}{U}\text{,}
\end{equation}
d'autre part on a
\begin{align}\label{eqdoro2}
q_2+\sqrt{1-q_1^2} \geq 1+\rho_0-d_0 && \text{ et } && q_2 \leq 1 +\rho_0 + d_0\text{,}
\end{align}
en sorte qu'en utilisant cela dans (\ref{eqdoro}) on obtient
\begin{equation}
  \label{eqdoro3}
  1+\rho_0-d_0 \leq e^4 (q_2-\hat q_2) \leq e^4(1+\rho_0+d_0 - \hat q_2)\text{.}
\end{equation}
On remplace $\hat q_2$ dans (\ref{eqdoro3}) par son expression déterminée par (\ref{eqdoro1})
ce qui nous donne l'inégalité
\begin{equation}
  \label{eqdoro5}
  1+\rho_0-d_0 \leq e^4\left(1+\rho_0+d_0 - (1+\rho_0)\frac{U -q_1}{U}\right)\text{,}
\end{equation}
en utilisant le fait que $q_1\leq d_0$, ceci nous donne l'inégalité (car $U\geq1$)
\begin{equation}
  \label{eqdoro6}
  1+\rho_0-d_0\leq e^4 d_0\left( 1 + \frac{1}{U}(1+\rho_0)\right) \leq e^4 d_0(2 +\rho_0)
\end{equation}
donc finalement
$$
d_0\geq \frac{1+\rho_0}{2e^4+1 +e^4\rho_0} \text{ ,} 
$$
qui est une fonction croissante en $\rho_0$, ce qui implique que $d_0\geq1/(2e^4+1)$.

En conclusion, on a toujours 
\begin{equation}
  \label{finetape1}
d_0\geq\frac{1}{2e^4+1}\text{.}  
\end{equation}

\textbf{Etape 2:} Pour tout $x \in \mathcal C$ et
pour toute droite $D$ issue de $x$, il existe un point $z\in D\cap\partial  \mathcal{C}$ à distance euclidienne
de $x$ majorée par 
$$
\Vert z-x\Vert_e \leq 5\sqrt{n}+3(2e^4+1)n \text{.}
$$

On commence par montrer qu le résultat est vrai pour le point 
$p$. Comme ce point est au centre de la boule $B_\mathcal{C}(p,1)$, toute droite issue de 
$p$ coupe la boule en un point $Y$ dont la distance à $p$ est exactement $1$.
Notons $A$ et $B$ les points d'intersection de cette droite avec 
$\partial \mathcal C$ et notons $y$, $a$ et $b$ les distances euclidiennes de $Y$, $A$ et $B$ à $p$.

On a l'égalité
\begin{equation}\label{distanceun}
\frac{1}{2}\ln \frac{y+a}{a} + \frac{1}{2}\ln\frac{b}{b-y} = 1\text{.}
\end{equation}
l'un des deux termes de l'égalité (\ref{distanceun}) est plus grand que $1/2$,
\begin{align*}
  \frac{1}{2}\ln \frac{y+a}{a}\geq \frac{1}{2}&& \text{ou}&& \frac{1}{2}\ln \frac{b}{b-y}\geq \frac{1}{2},
\end{align*}
et en utilisant le fait que $\n{Y}_e\leq \sqrt{n}$ on obtient
\begin{align}\label{equ2bruno}
  a\leq \frac{\sqrt{n}}{e-1}\leq 3\sqrt{n} && \text{ou} && b\leq\frac{e\sqrt{n}}{e-1}\leq3\sqrt{n}.
\end{align}

Considérons à présent un vecteur $v$ et la droite $D(x,v)$ passant par un point $x$ de $B_\mathcal{C}(p,1)$ 
et dirigée par $v$. Cette droite rencontre $\partial \mathcal{C}$ en deux point $A_x$ et $B_x$.

Le but est de contrôler la distance euclidienne entre $A_x$, $B_x$ et $x$. Pour cela,
on considère la droite issue de $p$ et parallèle à $D(x,v)$. 
On sait que sur cette droite il existe un point $A_p$, $A_p\in \partial \mathcal{C}$,
avec $\Vert p-A_p\Vert_e\leq 3\sqrt{n}$, ceci provient de l'inégalité (\ref{equ2bruno}).

Sans perte de généralité nous supposerons que $A_x$ est le point tel que 
$$ (v\cdot\vec{pA_p})(v\cdot\vec{xA_x})>0\text{.}$$

Si $p \in D(x,v)$,
la conclusion est immédiate, sinon on se place dans 
le plan engendré par ces deux droites parallèles.

Soit la perpendiculaire à $D(x,v)$ passant par $p$. Elle rencontre le bord
en deux points $D_p$ et $C_p$ tels que leur distance euclidienne à $p$ soit plus grande que $d_0$,
lui même minoré par $C_2=1/(2e^4+1)$ suivant l'étape 1. Sans perte de généralité nous noterons
$D_p$ le point tel que $\vec(pD_p)\cdot\vec(pA_x)<0$.

\begin{figure}[hbtp]
  \centering
  \includegraphics[scale=.6]{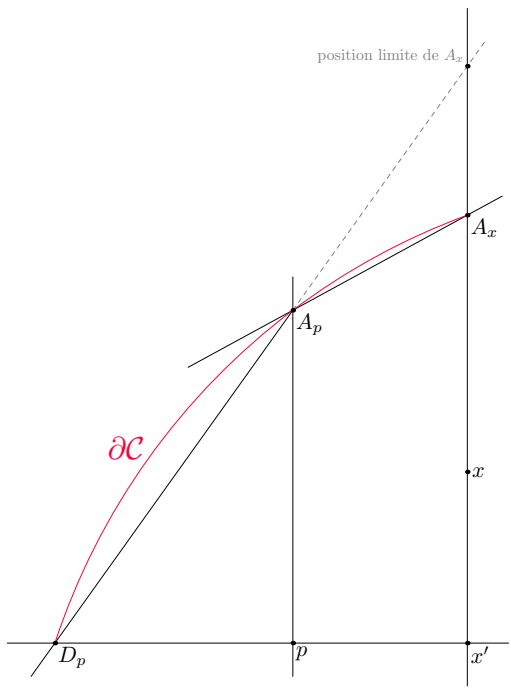}
 \caption{Position limite de $A_x$} \label{figtrois}
\end{figure}

Puisque $\partial \mathcal{C}$ passe par $A_p$ et $A_x$, par convexité, au dela du segment
$[A_p,A_x]$, $\partial \mathcal{C}$ se trouve en dessous de la droite passant pas $A_p$ et $A_x$.
De même en considérant le segment $[A_p,D_p]$. En sorte que  $A_x$ ne peut pas être plus haut
que la position critique où les droites $(A_pD_p)$ et $(A_xA_p)$ se confondent
(voir  figure \ref{figtrois}).

Plaçons nous donc dans cette position critique pour majorer la distance $xA_x$.
Soit $x'$ l'intersection de la droite $(xA_x)$ et de la droite $(pD_p)$
En comparant les longueurs des triangles homothétiques $(D_ppA_p)$ et $(D_px'A_x)$ on obtient

\begin{equation}
  \label{eqthales}
  \n{x'-A_x}_e=\n{p-A_p}_e \frac{\n{x'-D_p}_e}{\n{p-D_p}_e}
\end{equation}

Rappelons nous à présent que
\begin{enumerate}
\item $\n{p-A_p}_e \leq 3\sqrt{n}$ par (\ref{equ2bruno});
\item $\n{p-D_p}_e \geq \dfrac{1}{2e^4+1}$ par l'étape 1;
\end{enumerate}

ce qui donne une fois appliqué à l'équation (\ref{eqthales})
\begin{equation}
  \label{eqthales2}
  \n{x'-A_x}_e\leq 3\sqrt{n} \bigl(1+(2e^4+1)\n{x'-p}_e\bigr)
\end{equation}

Remarquons à présent que $x'$ est sur le cercle de diamètre $[x,p]$ dont la
largeur est majoré par $2\sqrt{n}$, puisque $x$ et $p$ sont dans le disque de
rayon $\sqrt{n}$ centré en l'origine. On en déduit donc que
\begin{align}
  \n{x'-p}_e\leq 2\sqrt{n} && \text{et} && \n{x-x'}_e\leq 2\sqrt{n}
\end{align}
ainsi en utilisant l'inégalité triangulaire dans l'espace euclidien et l'inégalité
(\ref{eqthales2}) on obtient finalement
\begin{equation}
  \label{eqfinetape2}
  \n{x-A_x}_e \leq \n_e{x-x'}+\n{x'-A_x}_e \leq 5\sqrt{n}+3(2e^4+1)n \text{.}
\end{equation}
ce qui est un contrôle de la distance de $x$ à $A_x$ en fonction de $n$.

\textbf{Etape 3:} Conclusion. On peut maintenant comparer la norme de Finsler et la norme
euclidienne en chaque point, comme annoncé au début de la preuve.

\smallskip
Pour majorer $\Vert v \Vert_{\mathcal C}$, $v$ vecteur tangent au point $x$,
 on utilise l'étape 1, qui nous donne via l'inégalité (\ref{finetape1})
\begin{align*}
  \n{x-x_+}_e \geq 2e^4+1 && \text{et} && \n{x-x_-}_e\geq 2e^4+1 \text{,}
\end{align*}
on en déduit donc
\begin{equation}
  \label{bilip1}
\Vert v \Vert_{\mathcal C}\leq(2e^4+1)\n{v}_e\text{.}  
\end{equation}

Pour minorer $\Vert v \Vert_{\mathcal C}$ on se place dans le (un des) plan(s) engendré(s) par $\vec{px}$ et $v$,
on travaille en dimension $2$ et on applique l'étape 2 
qui nous dit que soit $\n{x-x_-}_e$, soit $\n{x-x_+}_e$ est plus petit que $c(n)$.
On en déduit donc que
\begin{equation}
  \label{bilip2}
  \frac{1}{ 2c(n)}\n{v}_e\leq \Vert v \Vert_{\mathcal C}\text{.}
\end{equation}

Des inégalités (\ref{bilip1}) et (\ref{bilip2}), par intégration le long des chemins,
on en déduit que $(B_\mathcal{C}(p,1),\Vert . \Vert_{\mathcal C})$ est bilipschitz équivalente à $(B_\mathcal{C}(p,1),{\rm Eucl})$,
où la norme euclidienne que l'on a choisie est celle donnée par l'ellipsoïde de John associée à la boule unité. Une autre manière d'exprimer cela est de considérer l'application affine qui 
envoie l'ellipsoïde de John associé à la boule $B_\mathcal{C}(p,1)$ sur la boule unité de l'espace
euclidien usuel; c'est elle l'application bilipschitz.
Ceci permet de conclure notre preuve.
\end{proof}

\section{La constante de Cheeger et la $\delta$-hyperbolicité}

Nous commençons par démontrer le 

\begin{theo}\label{deltaIvol}
Un convexe du plan $\mathcal{C}$ qui, muni de sa métrique de Hilbert 
$d_\mathcal{C}$ est Gromov-hyperbolique, 
a sa constante de Cheeger strictement positive. Plus précisément : 
pour tout $\delta >0$, il existe $I>0$
tel que, pour tout convexe du plan $C$:
$$
(\mathcal{C},d_\mathcal{C}) \text{ est } \delta\text{-hyperbolique} \Rightarrow I^\infty(\mathcal{C})\geq I\text{.}
$$
\end{theo}

Pour démontrer ce théorème on va utiliser un résultat dû à
Jianguo~Cao~\cite{jcao}.

Nous avons au préalable besoin d'introduire une notion. Dans ce qui suit, $X$ est un espace métrique
de longueur.

\begin{defi}[quasi-pôle]
  Soit $\Omega\subset X$ un sous ensemble compact de $X$. Nous noterons $\mathcal{R}_\Omega$ l'union de tous les
rayons géodesiques minimisants émanant de $\Omega$. S'il existe une constante $C>0$ telle que le $C$-voisinage
de $\mathcal{R}_\Omega$ est égal à $X$, alors on dira que $X$ possède un quasi-pôle dans $\Omega$.
\end{defi}

Dans le cas des géométries de Hilbert, la situation est très simple: comme chaque
segment de droite est une géodésique minimisante, chaque point est un pôle: l'espace
$\mathcal C$ est réunion des rayons minimisants issus de tout point donné.

 \medskip
Nous pouvons énoncer alors le résultat important suivant qui s'appliquera dans notre cas

\begin{theo}[Main Theorem dans \cite{jcao}]
Soit $X$ une variété riemannienne non compacte (ou un graphe) à géométrie locale bornée admettant un quasi-pôle.
Si de plus $X$ est hyperbolique au sens de Gromov et le diamètre des composantes connexes du bord
à l'infini $X(\infty)$ est minoré par une constante strictement positive (respectivement à
une métrique de Gromov fixe sur le bord), alors $X$ a une constante de Cheeger strictement
positive.
\end{theo}

Remarquons que le concept de géométrie bornée qu'utilise Cao correspond   
  à la définition \ref {geobornee}. En particulier, sa preuve, donnée
dans le cas riemannien, s'adapte mot pour mot au cas finslérien qui nous concerne
ici. Cependant, on peut le voir directement via l'utilisation de l'ellipsoïde
de John: Soit $(\mathcal{C},d_\mathcal{C})$ un convexe du plan muni de sa métrique de Hilbert.
En tout point $p$ de $\mathcal{C}$, sur $T_p\mathcal{C}$ l'espace tangent en $p$,
on considère $J_p$ l'ellipsoïde de John associé à la boule unité de $\Vert . \Vert_{\mathcal C,p}$.
Toutes ces boules déterminent une nouvelle métrique, qui est riemannienne, et que l'on note $g_\mathcal{C}$.
Elle est bilipschitz-équivalente à la métrique de Finsler, puisque en tout point $p$ et $v\in T_p\mathcal{C}$ on a
$$
\frac{1}{\sqrt{n}}g_{\mathcal{C},p}(v,v)\leq \Vert v \Vert_{\mathcal C,p}^2 \leq g_{\mathcal{C},p}(v,v)\text{.}
$$

Cette métrique est en générale une métrique riemannienne $\mathcal{C}^0$. Pour pouvoir
appliquer ce résultat dans notre contexte, on utilise les trois remarques suivantes

\begin{itemize}
\item la possession d'un quasi-pôle
      pour une métrique de Finsler ou la propriété d'être à géométrie locale bornée sont
      des invariants bilipschitz.
\item la positivité de la constante de Cheeger pour la métrique riemannienne de John $g_\mathcal{C}$
      est équivalente à la positivité de la constante de Cheeger pour la métrique de Finsler initiale
      $\Vert . \Vert_{\mathcal C}$. 
\item Les bords à l'infini sont homéomorphes pour les deux métriques.
\end{itemize}

De plus, si la métrique de Hilbert est $\delta$-hyperbolique, le bord possède une seule composante
connexe de diamètre non nulle.

\medskip
Ainsi la seule chose non élémentaire à démontrer pour pouvoir appliquer
le théorème de Cao est le fait que les géométries de Hilbert sont à géométrie locale
bornée, ce qui faisait l'objet du théorème \ref{theoprincipal}. On obtient ainsi comme
corollaire le théorème \ref{equiv2}. La preuve du théorème \ref{equiv1}
se déduit alors directement de l'inégalité de Cheeger démontrée dans
\cite{cv1}, théorème 14.  

\section{La non nullité du bas du spectre n'implique pas la $\delta$-hyperbolicité}

Dans cette partie nous allons construire un convexe dans $(\R^3,\n{\cdot}_e)$, 
l'espace euclidien de dimension trois, dont la géométrie
de Hilbert n'est pas $\delta$-hyperbolique, mais dont le bas du spectre et la constante
de Cheeger ne sont cependant pas nuls. Intuitivement, en géométrie riemannienne,
il faut penser à une variété de Cartan-Hadamard à courbure presque partout
strictement négative, mais avec un plan euclidien qui fait que l'on est pas
Gromov-hyperbolique, tout en ayant un bas du spectre strictement positif.

\begin{prop}
Soit $(\mathcal{C},d_\mathcal{C})$ le cylindre de hauteur $2$ et de base le disque unité.
Sa constante de Cheeger et son bas du spectre sont non nuls. 
Cependant, puisqu'il n'est pas strictement convexe, il n'est pas hyperbolique au sens de Gromov.
\end{prop}

La remarque fondamentale est que lorsque l'on coupe le convexe $\mathcal{C}$ par des plans orthogonaux
 à l'axe de rotation, on obtient
un modèle de l'espace hyperbolique, pour lequel le bas du spectre est non nul (il vaut $1/4$). On va utiliser
cette propriété pour minorer le bas du spectre de $\mathcal{C}$ par une constante positive.

Pour cela, on aura besoin de comparer la forme volume de $\mathcal{C}$ et la forme volume
des sections orthogonales, obtenue
en restreignant la structure de Hilbert aux sections orthogonales.

Si on suppose que le cylindre est d'axe $Oz$, de base le disque unité, 
on note $D_t = \mathcal{C}\cap \{z=t\}$, l'intersection du cylindre $\mathcal{C}$ et du plan horizontal 
$\{z=t\}$, $-1 \leq t \leq 1$.

On paramétrise $\mathcal C$ de manière naturelle par les couples $p=(q,t)$,
où $q$ désigne un point du disque unité et $-1<t<1$. Le but est de 
comparer la forme volume de $\mathcal C$ en un point $p(q,t)$ avec la forme volume de
$D_t$ au point $q$. Pour cela, on doit comparer le volume euclidien des boules tangentes
$TB_{\mathcal C}(p,1)$ et $TB_{D_t}(q,1)$.

\begin{lemm}\label{cylindre} Pour tout $t\in \mathopen]-1,1\mathclose[$, notons $\alpha(t)=(1+t)(1-t)$. Alors
il existe deux constantes $C_1$ et $C_2>0$ telles que pour tout $p(q,t)$ on a
$$C_1 \alpha(t) \vole (TB_{D_t}(q,1)) \leq
\vole (TB_{\mathcal C}(p,1)) \leq C_2 \alpha(t)\vole (TB_{D_t}(q,1))\text{.}$$
\end{lemm}

Avant de démontrer ce lemme, montrons comment il permet de minorer le bas du spectre du cylindre.

\smallskip
\noindent\textbf{Minoration du spectre du cylindre $\mathcal C$.} 

Pour toute fonction
différentiable $f$ à support compact, la définition de la norme implique immédiatement
$$\Vert df\Vert_{\mathcal C}^* \geq \Vert df\Vert_{D_t}^* .$$

On a 
\begin{multline}
  \int_{\mathcal C}(\Vert df\Vert_{\mathcal C}^*)^2d\mu_{\mathcal C} =
\int_{\mathcal C}(\Vert df\Vert_{\mathcal C}^*)^2 \frac{\omega_3}{\vole (TB_{\mathcal C}(p,1))}d\vole(p)\geq\\
\int_{-1}^{1}dt\int_{D_t}(\Vert df\Vert_{D_t}^*)^2 
\frac{\omega_3}{\omega_2 C_2}
\frac{1}{\alpha(t)}
\frac{1}{\vole (TB_{D_t}(q,1))}d\vole(q)\text{;}
\end{multline}
En utilisant le fait que le spectre du plan hyperbolique $D_t$ est minoré par $1/4$, on obtient
\begin{equation}
\int_{\mathcal C}(\Vert df\Vert_{\mathcal C}^*)^2d\mu_{\mathcal C} \geq
\frac{\omega_3}{4\omega_2 C_2}
\int_{-1}^{1}\frac{dt}{\alpha(t)}\int_{D_t} f^2 \frac{1}{\vole (TB_{D_t}(q,1))}d\vole(q)\text{.}
\end{equation}
On a également
\begin{multline}
\int_{\mathcal C} f^2 d\mu_{\mathcal C} =
\int_{\mathcal C}f^2 \frac{\omega_3}{\vole (TB_{\mathcal C}(p,1))}d\vole(p)
\leq\\
\frac{\omega_3}{\omega_2 C_1}\int_{-1}^{1}\frac{dt}{\alpha(t)}\int_{D_t} f^2 
 \frac{1}{\vole (TB_{D_t}(q,1))}d\vole(q)\text{,}
\end{multline}
ce qui donne immédiatement une borne inférieure sur le spectre de $\mathcal C$.

 \bigskip
 \noindent
 \textbf{Preuve du lemme \ref{cylindre}.} En un point $p=(q,t)$, on considère les
 boules tangentes unités. Ce sont des convexes symétriques 
 $TB_{\mathcal C}(p,1)$ et $TB_{D_t}(q,1)$ que l'on supposera
 centrés à l'origine et qui vérifient 
 $TB_{D_t}(q,1) \subset TB_{\mathcal C}(p,1)$. Par ailleurs, comme la restriction de
 la métrique à $D_t$ est riemannienne, $TB_{D_t}(q,1)$ est une ellipse que
 l'on supposera située dans le plan $Oxy$.

 L'idée est de montrer que $TB_{\mathcal C}(p,1)$ est contenu dans un cylindre
de base un homothétique de $TB_{D_t}(q,1)$ et contient deux cônes de base
 $TB_{D_t}(q,1)$. Les hauteurs de ces deux objets étant controlées
par les faits suivant:

 \noindent
 \textbf{Fait 1 :}
 Si $\alpha(t)=(1+t)(1-t)$ les points $(0,0,\alpha(t))$ et $(0,0,-\alpha(t))$ sont
 dans $TB_{\mathcal C}(p,1)$. Notons que $\alpha(t)$ ne dépend pas de $q$.

 La remarque cruciale est alors 
 
 \smallskip
 \noindent
 \textbf{Fait 2 :} Au voisinage de ces deux points, le convexe 
 $TB_{\mathcal C}(p,1)$ est de la forme ${\{z=\alpha(t)\}}$ et ${\{z=-\alpha(t)\}}$.
 
 Cela se voit pas un calcul simple que nous expliciterons à la fin de la preuve.
 
 \medskip
 La conséquence est que, par convexité, la boule unité 
 $TB_{\mathcal C}(p,1)$ est située entre les plans $\{z=\alpha(t)\}$ et $\{z=-\alpha(t)\}$.
 
 Considérons les images, de $TB_{D_t}(q,1)$ sur $\{z=\alpha(t)\}$
 par la projection de centre $(0,0,-\alpha(t))$ et 
 sur $\{z=-\alpha(t)\}$
 par la projection de centre $(0,0,\alpha(t))$. Elles déterminent une section
 d'un cône qui contient le convexe $TB_{\mathcal C}(p,1)$.
 
 Cela signifie que le convexe $TB_{\mathcal C}(p,1)$ est contenu dans le produit
 $$
2TB_{D_t}(q,1)) × [-\alpha(t),\alpha(t)]\text{,}
$$
 et on en déduit que
 $$\vole (TB_{\mathcal C}(p,1)) \leq 8 \alpha(t)\vole(TB_{D_t}(q,1)).$$
 
 \medskip
 Par ailleurs, par convexité, les deux cônes de sommets 
 $(0,0,\alpha(t))$ et $(0,0,-\alpha(t))$ et de base $TB_{D_t}(q,1)$ sont contenus dans
 $TB_{\mathcal C}(p,1)$. On en déduit que
 $$\vole (TB_{\mathcal C}(p,1)) \geq 2/3 \vole (TB_{D_t}(q,1)).$$
 
 \medskip
 \noindent
 \textbf{Preuve du fait 2 :} C'est un fait général. On considère un point 
 $p$ d'un convexe $\mathcal C$ de l'espace $\R^3$. Sans limiter la généralité, on suppose
 $p$ situé à l'origine, et on suppose qu'au voisinage des points $(0,0,l_1)$
 et $(0,0,-l_2)$ avec $l_1$,$l_2>0$, le convexe $\mathcal C$ a pour bord les plans
 $\{z=l_1\}$ et $\{z=-l_2\}$. 
 
 Ainsi, les droites voisines de l'axe $Oz$, passant par l'origine, et faisant un 
 angle $\theta$ avec l'axe $Oz$, coupent ces deux plans à des distances respectivement égales
 à $\frac{l_1}{\cos \theta}$ et $\frac{l_2}{\cos \theta}$ de $p$. 
 
 Si $v$ est un vecteur directeur de la droite, sa norme de Finsler est
 
 $$\n{v}_{\mathcal C} = \frac{1}{2} \n{v}_{e} 
\biggl(\frac{\cos \theta}{l_1} + 
\frac{\cos \theta}{l_2}\biggr) {\rm }$$

et pour avoir $\n{v}_{\mathcal C}=1$, il faut que
$$\n{v}_{e} =\frac{1}{\cos \theta}\frac{2l_1l_2}{l_1+l_2}.$$
 
On voit que le lieu de ces points est bien le plan 
$z=\frac{2l_1l_2}{l_1+l_2}$ ou $z=-\frac{2l_1l_2}{l_1+l_2}$.
 
 \smallskip
\noindent\textbf{Minoration de la constante de Cheeger du cylindre $\mathcal C$.}

Remarquons que la méthode précédente permet également de montrer que la constante
de Sobolev
\begin{equation}
S^\infty(\mathcal{C}) =\inf \frac{\displaystyle{\int_\mathcal{C}{\n{df_p}_\mathcal{C}^*}~d\mu_\mathcal{C}(p)}}
{\displaystyle{\int_\mathcal{C} |f|(p)d\mu_\mathcal{C}(p)}}
\end{equation}
où l'infimum est pris sur toutes les fonctions 
lipschitziennes, non nulles, à support compact dans $\mathcal{C}$,
est également non nulle si $\mathcal{C}$ est un cylindre.

Cependant l'égalité entre la constante de Cheeger et la constante de Sobolev
n'a été démontrée que pour un convexe dont le bord est de régularité $C^1$. Il nous
faut donc justifier que cette égalité est encore valable dans notre cas.

Cette égalité provient de l'existence d'une formule de co-aire (voir théorème 28 dans
\cite{cv1} par exemple). C'est cette formule
qui dans l'article \cite{cv1} n'a été démontré que pour des bords $C^1$ et strictement convexes. 
Il nous faut donc justifier qu'une telle formule existe pour les convexes ayant un bord de régularité plus faible.

Rappelons que la formule de la co-aire  de \cite{cv1}, Théorème 25, nous dit que
si le bord du convexe $\mathcal{G}$ est $C^1$ et strictement convexe, 
alors pour toute fonction  $C^1$ à support compact $f$
l'égalité suivante est vérifiée
$$
\int_\mathcal{G} \n{df}^* d\mu_\mathcal{G} = \int_\R \int_{f^{-1}(t)} d\nu dt
$$
où $\nu$ est une mesure associée à $\mu_\mathcal{G}$ suivant la définition 22 dans \cite{cv1}
(c'est cette version faible qui est nécessaire dans la démonstration du théorème 28 de \cite{cv1}).

L'idée est d'approcher le cylindre par des domaines convexes de régularité $C^1$
et strictement convexes. Pour ces
domaines la formule est vraie. Il reste alors à vérifier que la formule passe à la limite.
Puisque $f$ est à support compact et vu la définition des mesures (Hilbert ou Holmes-Thompson), 
il suffit de vérifier que les normes Finsler convergent
uniformément sur le support de $f$ lorsque les convexes convergent vers le cylindre.
Le lemme suivant précise cela.
\begin{lemm}
  Soit $A$ un domaine compact de $\R^n$. Notons $\mathcal{O}(A)$ l'ensemble des
ouverts convexes bornés de $\R^n$ contenant $A$. A tout ouvert convexe borné $\mathcal{C}$ 
de $\mathcal{O}(A)$ on peut associer la norme de Hilbert induite par $\mathcal{C}$ sur
$A$, notons la $F_{\mathcal{C},A}$. Soit à présent une suite convergente $(\mathcal{C}_n)_{n\in\N}$ 
de convexes dans $\mathcal{O}(A)$ de limite $\mathcal{C}\in \mathcal{O}(A)$. Alors
la suite des applications de $A×\mathbb{S}^{n-1}$ dans $\R$ déterminées par
$$
(p,v)\mapsto \dfrac{F_{\mathcal{C}_n,A}(p,v)}{F_{\mathcal{C},A}(p,v)}
$$
converge uniformément vers $1$ sur $A×\mathbb{S}^{n-1}$.
\end{lemm}
\begin{proof}[Démonstration]
On va démontrer le résultat dans le cas où la suite de convexes
est décroissante au sens où $\mathcal{C}_{n+1}\subset\mathcal{C}_n$.

Considérons donc un point $p\in A \subset\mathcal{C}_n$ et un vecteur $v\in T_p\mathcal{C}$.
Soit $M_{n,A}\colon A×\mathbb{S}^1\to \R$ définie par 
$$
M_{n,A}(p,v)=\dfrac{F_{\mathcal{C}_n,A}(p,v)}{F_{\mathcal{C},A}(p,v)} \text{.}
$$
Si $m\geq n$, de l'inclusion $\mathcal{C}\subset\mathcal{C}_m\subset \mathcal{C}_n$, 
on obtient
$$
1\geq M_{m,A}(p,v)\geq M_{n,A} (p,v)\text{.}
$$
De plus $M_{\infty,A}(p,v)\equiv1$. 
La suite de fonctions  $(M_{n,A})_{n\in\N}$ converge simplement en croissant
vers $1$
(qui est une fonction continue)  sur l'ensemble compact $A× \mathbb{S}^1$. 
On en déduit qu'en réalité la convergence est uniforme sur $A× \mathbb{S}^1$.
\end{proof}
 
\begin{coro}
  Sous les mêmes hypothèses que le lemme précédant, en notant $h_n$ et $h_\infty$
les densités par rapport à la mesure de Lebesgue 
des mesures de Hilbert associés repectivement à $\mathcal{C}_n$
et $\mathcal{C}$, alors $h_n$ converge uniformément vers $h_\infty$ sur $A$. 
\end{coro}
\begin{proof}[Démonstration]
   soit
$$
K(n)=\sup\{M_{n,A}(p,v) \mid (p,v) \in A×\mathbb{S}^{n-1}\}\text{.}
$$
Par continuité de $M_{n,A}$ 
et compacité de  $A×\mathbb{S}^1$, $K(n)$ est fini.
Le lemme précédant implique que $K(n)$ converge en croissant vers $1$.
On en déduit que pour tout $\varepsilon>0$, à partir d'un certain rang
$$
(1+\varepsilon)^{-1/n}\leq K(n)<1
$$
et donc pour tout $p \in A$
$$
h_\infty(p)(1+\varepsilon) \geq h_n(p)>h_\infty(p)\text{,}
$$
ce qui implique bien la convergence uniforme de $h_n$ sur $A$.
\end{proof}
